\documentclass[11pt,reqno]{amsart}
\usepackage{graphicx}
\usepackage{amsmath}
\usepackage{eucal}
\usepackage{amscd}
\usepackage{amssymb}
\usepackage{latexsym}
\usepackage{amsfonts}
\newtheorem{thm}{Theorem}[section]
\newtheorem{lemma}{Lemma}[section]
\newtheorem{prop}{Proposition}[section]
\numberwithin{equation}{section}

\title{Standard 2D Crystalline Patterns and Rational Points in Complex Quadrics}

\author{Toshikazu Sunada}
\address{
Department of Mathematics, Meiji University, 
Higashimita 1-1-1, Tama-ku, Kawasaki, 214-8571 Japan
}
\email{sunada@isc.meiji.ac.jp}
\date{}
\begin{document}
\maketitle
\pagestyle{myheadings}

\markright{Standard 2D Crystalline Patterns}
\begin{abstract}
A certain Diophantine problem and 2D crystallography are linked through the notion of {\it standard realizations} which was introduced originally in the study of random walks.
In the discussion, a complex projective quadric defined over $\mathbb{Q}$ is associated with a finite graph. ``Rational points" on this quadric turns out to be related to standard realizations of 2D crystal structures. In the last section, it is observed that the number of rational points which correspond to periodic tilings is finite.
\end{abstract}

\section{Introduction}
The standard crystalline pattern is a synonym of the {\it standard realization} (or canonical placement) of a crystal structure introduced in \cite{sk2} which  gives the most symmetric crystalline shape among all possible realizations, and is characterized uniquely (up to similar transformations) as a minimizer of a certain energy functional. 

\begin{figure}[htbp]
\begin{center}
\includegraphics[width=.6\linewidth]
{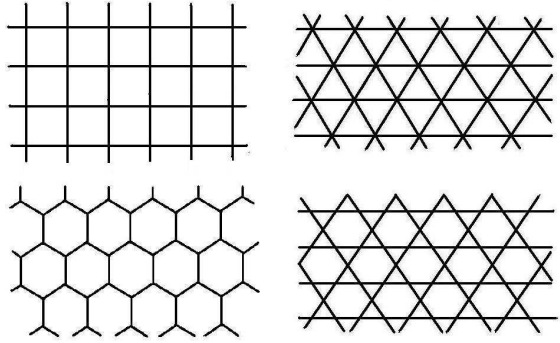}
\end{center}
\caption{Typical 2D standard realizations}\label{fig:classical}
\end{figure}

Figure \ref{fig:classical} exhibits some of 2D examples which are familiar to scientists; the square lattice, honeycomb (regular hexagonal lattice), equilateral triangular lattice, and regular kagome lattice. In the 3D case, Diamond, Lonsdaleite\footnote{This carbon allotrope, formed when meteorites containing graphite strike the earth, is named in honor of crystallographer Kathleen Lonsdale, also referred to as the {\it hexagonal diamond}.}, and the nets associated with the face-centered and body-centered lattices 
give standard crystalline patterns. This notion, via the elementary theory of homology and covering spaces, offers an effective method in a systematic enumeration and design of crystal structures (see \cite{su4})\footnote{Crystallographers proposed a similar idea in \cite{d-1}, \cite{d1}, \cite{eon}, \cite{eon1}}. Remarkably, it shows up in asymptotic behaviors of simple random walks (\cite{sk1}), and is also related to a combinatorial analogue of Abel-Jacobi maps (\cite{su5}).

The aim of this paper is to observe another interesting aspect of standard realizations; specifically we establish a link between standard 2-dimensional crystalline patterns and ``rational points" of certain complex quadrics defined over $\mathbb{Q}$. A rational point we mean here is a point in a complex projective space each of whose homogeneous coordinate is represented by a number in an {\it imaginary quadratic field}.   
For instance, the regular kagome lattice (the lower right in Fig.~\!\ref{fig:classical}) corresponds to the $\mathbb{Q}(\sqrt{-3})$-rational point 
$$
\Big[\frac{1+\sqrt{-3})}{2}, \frac{1-\sqrt{-3}}{2}, -1, \frac{1+\sqrt{-3}}{2}, \frac{1-\sqrt{-3}}{2}, -1\Big]
$$
(or its complex conjugate) 
of the 2-dimensional projective quadric defined over $\mathbb{Q}$
\begin{eqnarray*}
&&\{[z_1,z_2,z_3, z_4,z_5,z_6]\in P^5(\mathbb{C});~z_1{}^2+\cdots+z_6{}^2=0,\\
&& \qquad z_1+z_6=z_3+z_4,~ z_2+z_4=z_1+z_5,~z_3+z_5=z_2+z_6\},
\end{eqnarray*}
which is biregular (over $\mathbb{Q}$) to
\begin{eqnarray*}
&& \{[w_1,w_2,w_3,w_4]\in P^3(\mathbb{C})|~3w_1{}^2+2w_2{}^2+2w_3{}^2+3w_4{}^2\\
&& \qquad \qquad -2w_1w_2-4w_1w_4-2w_1w_3
  +2w_2w_4+2w_3w_4=0\}.
\end{eqnarray*}

The way to give this correspondence is not {\it ad hoc}. Indeed there is a systematic procedure to associate a complex projective quadrics ${\bf Q}(X_0)$ defined over $\mathbb{Q}$ with a {\it finite graph} $X_0$, which is an interesting object for its own sake. In our setting, $X_0$ turns up as the quotient graph by the translational action of a lattice group on the crystal structure. The point corresponding to the standard crystalline pattern is obtained by taking the intersection of this quadric and a certain (projective) line. It is also observed that every rational point on the quadric yields a standard 2D crystalline pattern. 

Incidentally all the examples given in Fig.~\!\ref{fig:classical} happen to be derived from {\it tilings} of plane. As a matter of fact, it is rather rare to find tilings among standard 2D crystalline patterns; namely as will be proved in the last section, there are only finitely many rational points on ${\bf Q}(X_0)$ which correspond to tilings.

\section{Standard realizations}\label{sec:standard}
We first review some basic results about standard realizations (see \cite{su5} for the terminology used in the present article).

A {\it graph}\index{graph} is represented by an ordered pair $X = (V,E)$ of the set of {\it vertices} $V$ and the set of all {\it directed edges} $E$ (note that each edge has just two directions, which are to be expressed by arrows). For an  directed edge $e$, we denote by ${\it o}(e)$ the {\it origin}, and by ${\it t}(e)$ the {\it terminus}. The inversed edge of $e$ is denoted by $\overline{e}$. With these notations, we have $o(\overline{e})=t(e)$, $t(\overline{e})=o(e)$. We also use the notation $E_x$ for the set of directed edges $e$ with $o(e)=x$. Throughout, the degree ${\rm deg}~\!x=|E_x|$ is assumed to be greater than or equal to three for every vertex $x$.

The network associated with a $d$-dimensional crystal is identified with a periodic realization of a $d$-dimensional {\it topological crystal} (an infinite-fold covering graph) $X=(V,E)$ over a finite graph $X_0=(V_0,E_0)$ whose covering transformation group is a free abelian group $L$ of rank $d$ ($d=2$ or $3$ when we are handling a network associated with a real crystal). Theory of covering spaces tells us that there is a subgroup $H$ (called a {\it vanishing subgroup}) such that $H_1(X_0,\mathbb{Z})/H=L$ (note that $H$ is a direct summand of $H_1(X_0,\mathbb{Z})$). Actually the topological crystal $X$ is the quotient graph of the maximal abelian covering graph $X_0^{\rm ab}$ over $X_0$ modulo $H$. We call $X_0^{\rm ab}$ the {\it maximal topological crystal} over $X_0$.
We denote by $\mu:H_1(X_0,\mathbb{Z})\longrightarrow L$ the canonical homomorphism.

Precisely speaking, a periodic realization is a piecewise linear map $\varPhi:X\longrightarrow \mathbb{R}^d$ satisfying 
$$
\varPhi(\sigma x)=\varPhi(x)+\rho(\sigma)\qquad (\sigma \in L),
$$
where $\rho:L\longrightarrow \mathbb{R}^d$ is an injective homomorphism whose image is a lattice in $\mathbb{R}^d$.\footnote{The network constructed in this way could be ``degenerate" in the sense that different vertices of $X$ are realized as one points, or different edges overlap in $\mathbb{R}^d$. But we shall not exclude these possibilities.} We call $\rho$ (resp. $\rho(L)$) the {\it period homomorphism} (resp. the {\it period lattice}) for $\varPhi$.

By putting ${\bf v}(e)=\varPhi\big(t(e)\big)-\varPhi\big(o(e)\big)$ ~$(e\in E)$, we obtain a $L$-invariant function ${\bf v}$ on $E$ which we may identify with a 1-cochain ${\bf v}\in C^1(X,\mathbb{R}^d)$ with values in $\mathbb{R}^d$. 
Since ${\bf v}$ determines completely $\varPhi$ (up to parallel translations), we shall call ${\bf v}$ the {\it building cochain}\footnote{In \cite{su4}, \cite{su5}, the term ``building block" is used. The idea to describe crystal structures by using finite graphs together with vector labeling is due to \cite{chung}} of $\varPhi$. 
One can check that if we identify the cohomology class $[{\bf v}]\in H^1(X_0,\mathbb{R}^d)$ with a homomorphism of $H_1(X_0,\mathbb{Z})$ into $\mathbb{R}^d$ (the {\it duality of cohomology and homology}), then $[{\bf v}]=\rho\circ \mu$. In particular, ${\rm Ker}~\![{\bf v}]=H$ and ${\rm Image}~\![{\bf v}]=\rho(L)$. 

\begin{lemma}{\rm (\cite{su4})}
Giving a periodic realization of a topological crystal over $X_0$ is equivalent to giving a 1-cochain ${\bf v}\in C^1(X_0,\mathbb{R}^d)$ such that the image of the homomorphism $[{\bf v}]:H_1(X_0,\mathbb{Z})\longrightarrow \mathbb{R}^d$ is a lattice in $\mathbb{R}^d$. 
\end{lemma}

Among all periodic realizations of $X$, there is a ``standard" one which is characterized uniquely (up to similar transformations) by the following two conditions:

\smallskip

(1) ({\bf Harmonicity})
\begin{equation}\label{eq:harmonic}
\displaystyle\sum_{e\in E_{0x}}{\bf v}(e)={\bf 0} \quad (x\in V_0),
\end{equation}

\smallskip
(2) ({\bf Tight-frame condition}\footnote{The term ``tight frame" is the terminology in wavelet analysis.}) There exists a positive constant $c$ such that
\begin{equation}\label{eq:standard}
\sum_{e\in E_0} \langle{\bf x},{\bf v}(e) \rangle {\bf v}(e)=c{\bf x}\quad ({\bf x}\in \mathbb{R}^d).
\end{equation}

\medskip

In the coordinate form, (\ref{eq:standard}) is written as
$$
\sum_{e\in E_0}v_i(e)v_j(e)=c\delta_{ij},
$$
where ${\bf v}={}^t\big(v_1(e),\ldots,v_d(e)\big)$. In particular
\begin{equation}\label{eq:cd}
\sum_{e\in E_0}\|{\bf v}(e)\|^2=cd.
\end{equation}
(\ref{eq:standard}) is also equivalent to 
$$
\sum_{e\in E_0} \langle{\bf x},{\bf v}(e) \rangle^2=c\|{\bf x}\|^2.
$$

The constant $c$ does not play a role since we are considering similarity classes of realizations; later we handle the case $c=2$.

\begin{prop} {\rm (\cite{sk2}, \cite{sk11})}

{\rm (1)} The standard realization is the unique minimizer, up to similar transformations, of the energy\footnote{The energy defined here is similarity-invariant.} 
$$
\mathcal{E}(\varPhi)={\rm vol}\big(\mathbb{R}^d/\rho(L)\big)^{-2/d}\sum_{e\in E_0}\|{\bf v}(e)\|^2.
$$

{\rm (2)} Let $\varPhi:X\longrightarrow \mathbb{R}^d$ be the standard realization. Then there exists a homomorphism $\kappa$ of the automorphism group ${\rm Aut}(X)$ of $X$ into the congruence group $M(d)$ of $\mathbb{R}^d$ such that

{\rm (a)} when we write $\kappa(g)=\big(A(g),b(g)\big)\in O(d)\times \mathbb{R}^d$, we have
$$
\varPhi(gx)=A(g)\varPhi(x)+b(g) \qquad (x\in V), 
$$

{\rm (b)} the image $\kappa\big({\rm Aut(X)}\big)$ is a crystallographic group\index{crystallographic group}.
\end{prop}

\noindent {\bf Remark}~ Equation (\ref{eq:harmonic}) says that the cochain ${\bf v}$ is ``harmonic" in the sense that $\delta{\bf v}=0$ where $\delta:C^1(X_0,\mathbb{R}^d)\longrightarrow C^1(X_0,\mathbb{R}^d)$ is the adjoint of the coboundary operator $d:C^0(X_0,\mathbb{R}^d)\longrightarrow C^1(X_0,\mathbb{R}^d)$ with respect to the natural inner products in $C^i(X_0,\mathbb{R}^d)$. Using a discrete analogue of the Hodge--Kodaira theorem, one can prove that the correspondence ${\bf v}\in {\rm Ker}~\!\delta\mapsto [{\bf v}]\in H^1(X_0,\mathbb{R}^d)$ is a linear isomorphism (hence ${\rm dim}~\!{\rm Ker}~\!\delta=db_1(X_0)$ where $b_1(X_0)$ is the first Betti number of $X_0$). Thus given $\rho$, there is a unique harmonic cochain ${\bf v}$ with $[{\bf v}]=\rho\circ \mu$. A realization satisfying (\ref{eq:harmonic}) is said to be a {\it harmonic realization} \cite{sk2} (or an equilibrium placement \cite{d-1}), which is characterized as a minimizer of $\mathcal{E}$ when $\rho$ is fixed.

\section{Albanese tori}\label{sec:alb}
The building cochain ${\bf v}$ for the standard realization of $X$ with $c=2$ in Eq.~\!(\ref{eq:standard}) is explicitly constructed in the following way (\cite{su5}). 

First we provide $H_1(X_0,\mathbb{R})$ with a natural inner product (which allows us to identify $H_1(X_0,\mathbb{R})$ with the Euclidean space $\mathbb{R}^b$, $b=b_1(X_0)$). 
For this sake, we start with an inner product on $C_1(X_0,\mathbb{R})$. 

For $e,e'\in E_0$, we set
$$
\langle e, e'\rangle=\begin{cases}
1   & (e'=e)\\
-1 &  (e'=\overline{e})\\
0  &  (\text{otherwise})
\end{cases},
$$
which extends to an inner product on $C_1(X_0,\mathbb{R})$ in a natural manner.
Restricting this inner product to the subspace $H_1(X_0,\mathbb{R})~(={\rm Ker}~\big(\partial: C_1(X_0,$ $\mathbb{R})\longrightarrow C_0(X_0,$ $\mathbb{R})\big))$, we get an Euclidean structure on $H_1(X_0,\mathbb{R})$.

Let $P_{\rm ab}:C_1(X_0,\mathbb{R})\longrightarrow H_1(X_0,\mathbb{R})$ be the orthogonal projection, and put ${\bf v}_{\rm ab}(e)=P_{\rm ab}(e)$, regarding each edge as a $1$-chain. Note that $[{\bf v}^{\rm ab}]:H_1(X_0,\mathbb{Z})\longrightarrow H_1(X_0,\mathbb{R})$ coincides with the injection. One can check that ${\bf v}_{\rm ab}$ is the building cochain of the standard realization of the maximal topological crystal $X_0^{\rm ab}$ (with $c=2$). 

Let $X$ be a topological crystal over $X_0$ corresponding to a vanishing subgroup $H$ of $H_1(X_0,\mathbb{Z})$.
Let $H_{\mathbb{R}}$ be the subspace of $H_1(X_0,\mathbb{R})$ spanned by $H$, 
and $H_{\mathbb{R}}^{\perp}$ the orthogonal complement of $H_{\mathbb{R}}$ in $H_1(X_0,\mathbb{R})$:
$$
H_1(X_0,\mathbb{R})=H_{\mathbb{R}}\oplus H_{\mathbb{R}}^{\perp}.
$$
Then ${\rm dim}~\hspace{-0.05cm}H_{\mathbb{R}}^{\perp}={\rm rank}~\hspace{-0.05cm}L=d$. By choosing an orthonormal basis of $H_{\mathbb{R}}^{\perp}$, we identify $H_{\mathbb{R}}^{\perp}$ with the Euclidean space $\mathbb{R}^d$.

Let $P:H_1(X_0,\mathbb{R})\longrightarrow H_{\mathbb{R}}^{\perp}=\mathbb{R}^d$ be the orthogonal projection. If we put ${\bf v}(e)=P\big({\bf v}_{\rm ab}(e)\big)$, then 
we find that ${\bf v}$ gives the building cochain of the standard realization of $X$ (with $c=2$). As shown in \cite{su5}, one may reduce the construction of ${\bf v}$ to an elementary computation of matrices.

We now consider two flat tori 
\begin{eqnarray*}
&&A(X_0)=H_1(X_0,\mathbb{R})/H_1(X_0,\mathbb{Z}),\\
&&A(X_0,H)=\mathbb{R}^d/{\rm Image}~\![{\bf v}],
\end{eqnarray*}
which, in view of analogy with classical algebraic geometry, are called the {\it Albanese torus} of $X_0$ and the {\it generalized Albanese torus} of $(X_0,H)$, respectively. The projection $P$ induces the exact sequence
\begin{equation}\label{exactexact}
0\longrightarrow H_{\mathbb{R}}/H\longrightarrow A(X_0)\stackrel{p}{\longrightarrow}A(X_0,H)\longrightarrow 0,
\end{equation}

In the following proposition, $\kappa(X_0)$ denotes the {\it tree number} of $X_0$, the number of spanning trees in $X_0$.

\begin{prop} {\rm (1) (\cite{sk12})} 
$
{\rm vol}\big(A(X_0)\big)=\kappa(X_0)^{1/2}.
$

$(2)$ ${\rm vol}\big(A(X_0,H)\big)=\kappa(X_0)^{1/2}/{\rm vol}(H_{\mathbb{R}}/H).$

\end{prop}

The second claim is a consequence of the exact sequence (\ref{exactexact}).

The volume ${\rm vol}(H_{\mathbb{R}}/H)$ is computed as 
$$
{\rm vol}(H_{\mathbb{R}}/H)=\det (\langle \alpha_i, \alpha_j\rangle)^{1/2},
$$
where $\{\alpha_1,\ldots, \alpha_{b-d}\}$ is a $\mathbb{Z}$-basis of $H$. Putting $I(H)=\det (\alpha_i\cdot \alpha_j)$,  
we shall call $I(H)$ the {\it intersection determinant} for $H$, which is evidently a positive integer. We thus have the following formula.
$$
{\rm vol}\big(A(X_0,H)\big)=\kappa(X_0)^{1/2}I(H)^{-1/2},
$$
and hence, appealing to (\ref{eq:cd}), we obtain
\begin{equation}\label{eq:energymini}
\min_{\varPhi}~\!\mathcal{E}(\varPhi)=2d\kappa(X_0)^{-1/d}I(H)^{1/d},
\end{equation}
where $\varPhi$ runs over all periodic realizations of $X$.

\section{Complex quadrics associated with finite\\ graphs}
We now embark on a new enterprise. We confine ourselves to 2D standard realizations.

Let $X_0=(V_0,E_0)$ be a finite connected graph such that $b_1(X_0)$ is greater than or equal to 2.  Put
$$
\mathbb{H}=\{{\bf z}\in C^1(X_0,\mathbb{C})|~\sum_{e\in E_{0x}}{\bf z}(e)=0~ (x\in V_0)\}.
$$
This is nothing but the space of harmonic cochains (we identify $\mathbb{R}^2$ with $\mathbb{C}$), so we find 
$$
{\rm dim}_{\mathbb{C}}\mathbb{H}=b_1(X_0)
$$
(see Remark at the end of Sect.~\!\ref{sec:standard}). 
We denote by $\mathbb{P}(\mathbb{H})$ the projective space associated with $\mathbb{H}$; that is, $\mathbb{P}(\mathbb{H})$ is the orbit space $(\mathbb{H}\backslash\{{\bf 0}\})/\mathbb{C}^{\times}$ by the natural action of the multiplicative group $\mathbb{C}^{\times}=\mathbb{C}\backslash\{0\}$. For ${\bf z}\neq {\bf 0} (\in \mathbb{H})$, we use the notation $\langle{\bf z}\rangle\in \mathbb{P}(\mathbb{H})$ for the orbit containing ${\bf z}$. It should be noted that an orientation-preserving similar transformation in $\mathbb{R}^2$ is identified with the multiplication by a non-zero complex number.

\begin{lemma}
A cochain ${\bf z}\in C^1(X_0,\mathbb{C})=C^1(X_0,\mathbb{R}^2)$ satisfies the tight-frame condition if and only if $\displaystyle\sum_{e\in E_0}{\bf z}(e)^2=0$.
\end{lemma}

\noindent {\it Proof}~ Put ${\bf z}(e)=a(e)+b(e)\sqrt{-1}$, and ${\bf x}={}^t(x,y)$. Then 
$$
\sum_{e\in E_0}\langle {\bf z}(e),{\bf x}\rangle^2=
\Big(\sum_{e\in E_0}a(e)^2\Big)x^2+2\Big(\sum_{e\in E_0}a(e)b(e)\Big)xy+
\Big(\sum_{e\in E_0}b(e)^2\Big)y^2.
$$ 
Thus ${\bf z}$ satisfies the tight-frame condition if and only if 
$$
\sum_{e\in E_0}a(e)^2=\sum_{e\in E_0}b(e)^2, \quad \sum_{e\in E_0}a(e)b(e)=0.
$$
Since
$$
\sum_{e\in E_0}{\bf z}(e)^2=\sum_{e\in E_0}\big(a(e)^2-b(e)^2+2\sqrt{-1}a(e)b(e)\big),
$$
we get the claim.\hfill $\Box$

\medskip 

In view of this lemma, it is natural to consider the quadric defined by
$$
{\bf Q}(X_0)=\Big\{
\langle{\bf z}\rangle\in {\bf P}(\mathbb{H})\big|~\sum_{e\in E_0}{\bf z}(e)^2=0
\Big\}
$$

Let $H (\subset H_1(X_0,\mathbb{Z}))$ be a vanishing subgroup such that ${\rm rank}~\!H_1(X_0,\mathbb{Z})/H=2$. We consider the subspace of $\mathbb{H}$ defined by
$$
\mathbb{W}_H=\{{\bf z}\in \mathbb{H}|~[{\bf z}](\alpha)=0~ (\alpha\in H)\}
$$
(remember that the cohomology class $[{\bf z}]$ is identified with a homomorphism of $H_1(X_0,\mathbb{Z})$ into $\mathbb{C}$). 
Clearly ${\rm dim}_{\mathbb{C}}\mathbb{W}=2$ so that $\mathbb{L}_H={\bf P}(\mathbb{W}_H)$ is a (projective) line  in ${\bf P}(\mathbb{H})$.

Let ${\bf z}_H\in \mathbb{H}$ be the building cochain of the standard realization of $X$ corresponding to $H$ which was constructed in the previous section. We fix a complex structure on $H_{\mathbb{R}}^{\perp}$, and identify $H_{\mathbb{R}}^{\perp}$ with $\mathbb{C}$. It is evident that both ${\bf z}_H$ and its complex conjugate $\overline{{\bf z}_H}$ belong to $\mathbb{W}_H$. 

\begin{lemma}  $\mathbb{W}_H$ is spanned by ${\bf z}_H$ and $\overline{{\bf z}_H}$.
\end{lemma}

\noindent{\it Proof}~ It is enough to check that ${\bf z}_H$ and $\overline{{\bf z}_H}$ are linearly independent over $\mathbb{C}$. If $\overline{{\bf z}_H}=w{\bf z}_H$ for some $w\in \mathbb{C}$, then
$$
\sum_{e\in E_0}|{\bf z}_H(e)|^2=w\sum_{e\in E_0}{\bf z}_{H}(e)^2=0
$$
so that ${\bf z}_H={\bf 0}$, thereby a contradiction.
\hfill$\Box$

\medskip

The following theorem tells us that the standard realization of the 2D topological crystal over $X_0$ associated with the vanishing group is obtained as intersection of the quadric ${\bf Q}(X_0)$ and the line $\mathbb{L}_H$.

\begin{thm}
${\bf Q}(X_0)\cap \mathbb{L}_H=\{\langle {\bf z}_H\rangle, \langle \overline{{\bf z}_H}\rangle\}$. 
\end{thm}

\begin{figure}[htbp]
\begin{center}
\includegraphics[width=.37\linewidth]
{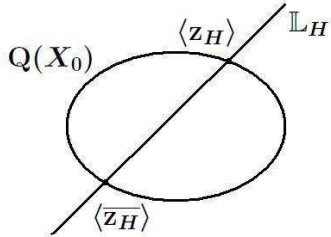}
\end{center}
\caption{Quadric and line}\label{fig:inter}
\end{figure}

\noindent{\it Proof}~ We only have to check that the line $\mathbb{L}_H$ is not contained in ${\bf Q}(X_0)$. But this follows immediately from $\sum_{e\in E_0}|{\bf z}_H(e)|^2\neq 0$ and 
$$
\sum_{e\in E_0}\big(a{\bf z}_H(e)+b\overline{{\bf z}_H(e)}\big)^2
=2ab\sum_{e\in E_0}|{\bf z}_H(e)|^2\quad (a,b\in\mathbb{C}).
$$
\hfill$\Box$

\medskip 

To give a coordinate form of ${\bf Q}(X_0)$, choose an orientation $E_0^o=(e_1,\ldots,e_N)$ of $X_0$, and put
\begin{eqnarray*}
&&E_{0x}^{\rm in}=\{e\in E_0^0|~t(e)=x\},\\
&&E_{0x}^{\rm out}=\{e\in E_0^0|~o(e)=x\}.
\end{eqnarray*}
Then ${\bf Q}(X_0)$ is identified with
\begin{eqnarray*}
&&\Big\{[z_1,\ldots,z_N]\in P^{N-1}(\mathbb{C})|~z_1{}^2+\cdots+z_N^{2}=0,\\
 &&\qquad \qquad \qquad
\sum_{i;e_i\in E_{0x}^{\rm in}}z_i=\sum_{j;e_j\in E_{0x}^{\rm out}}z_j ~~(x\in V_0)
\Big\}.
\end{eqnarray*}
via the correspondence ${\bf z}\mapsto (z_1,\ldots,z_N)$ given by $z_i={\bf z}(e_i)$. 
Note the equation
\begin{equation}\label{eq:kir}
\sum_{i;e_i\in E_{0x}^{\rm in}}z_i=\sum_{j;e_j\in E_{0x}^{\rm out}}z_j
\end{equation}
is equivalent to $\displaystyle\sum_{e\in E_{0x}}{\bf z}(e)=0$, and is a complex version of {\it Kirchhoff's law} in the theory of electric circuits stating that the amount of in-coming currents at a node is equal to the amount of out-going current. 

\begin{prop}
${\bf Q}(X_0)$ is (biregular over $\mathbb{Q}$ to) a non-singular quadric defined over $\mathbb{Q}$ of dimension $b_1(X_0)-2$.
\end{prop}
\noindent{\it Proof}~ Put 
$$
W=\Big\{{}^t(z_1,\ldots,z_N)|~\sum_{i;e_i\in E_{0x}^{\rm in}}z_i=\sum_{j;e_j\in E_{0x}^{\rm out}}z_j ~~(x\in V_0)\Big\},
$$
which is a coordinate form of $\mathbb{H}$.
One can find a $N\times b$ matrix $A=(a_{ij})$ with rational entries such that the correspondence ${}^t(w_1,\ldots,w_b)\mapsto {}^t(z_1, \ldots,z_N)$ given by $z_i=\sum_{j=1}^b a_{ij}w_j$ is a linear isomorphism of $\mathbb{C}^b$ onto $W$. Then 
$$
F(w_1,\ldots,w_b)=\sum_{i=1}^N\Big(\sum_{j=1}^ba_{ij}w_j\Big)^2
$$
is obviously a positive definite quadratic form. Thus the quadric 
$$
Q=\{[w_1,\ldots,w_b]\in P^{b-1}(\mathbb{C})|~F(w_1,\ldots,w_b)=0\}
$$
is non-singular, and is biregular (over $\mathbb{Q}$) to ${\bf Q}(X_0)$. 
\hfill $\Box$

\medskip

To describe $\{\langle{\bf z}_H\rangle, \langle\overline{{\bf z}_H}\rangle\}$ in the coordinate form, let $\alpha_1,\ldots,\alpha_{b-2}$ be a free $\mathbb{Z}$-basis of the vanishing subgroup $H$.
Write 
$$
\alpha_i=\sum_{j=1}^Na_{ij}e_j \quad (i=1,\ldots,b-2).
$$
with $a_{ij}\in \mathbb{Z}$. Then
$$
\mathbb{W}_H=\Big\{{}^t(z_1,\ldots,z_N)\big|~\sum_{j=1}^Na_{ij}z_j=0 ~ (i=1,\ldots,b-2)
\Big\}.
$$
We thus have

\begin{thm} 
\begin{eqnarray*}
&&\{\langle {\bf z}_H\rangle, \langle\overline{{\bf z}_H}\rangle\}=\Big\{[z_1,\ldots,z_N]\in P^{N-1}(\mathbb{C})|~z_1{}^2+\cdots+z_N^{2}=0,\\
 &&\qquad \qquad \qquad\quad
\sum_{i;e_i\in E_{0x}^{\rm in}}z_i=\sum_{j;e_j\in E_{0x}^{\rm out}}z_j ~~(x\in V_0),\\
&&\qquad \qquad \qquad\qquad \sum_{j=1}^Na_{ij}z_j=0 ~ (i=1,\ldots,b-2)
\Big\}.
\end{eqnarray*}
\end{thm}




\section{Rational points}
From the previous theorem, it follows that ${\bf z}_H$ and $\overline{{\bf z}_H}$ are obtained by solving a quadratic equation of the form $az_1{}^2+bz_1z_2+cz_2{}^2=0$ ~($a,b,c\in \mathbb{Q}$). Hence we conclude that there exists a positive square free integer $D$ such that $\{\langle {\bf z}_H\rangle, \langle\overline{{\bf z}_H}\rangle\}\subset P^{N-1}\big(\mathbb{Q}(\sqrt{-D})\big)$. 
The following lemma will give another proof for this fact.

\begin{lemma}
Suppose that $\{z_1,\ldots,z_N\}$ satisfies $z_1{}^2+\cdots+z_N^{2}=0$. Then $\{z_1,\ldots,z_N\}$ generates a lattice in $\mathbb{C}$ if and only if  there exists a positive square free integer $D$ such that $[z_1,\ldots,z_N] \in P^{N-1}\big(\mathbb{Q}(\sqrt{-D})\big)$. 
\end{lemma}

\noindent{\it Proof}~ The proof is fairly elementary. Suppose that $\{z_1,\ldots,z_N\}$ generates a lattice
whose $\mathbb{Z}$-basis $w_1,w_2$. Then there exist integers $a_i,b_i$ such that $z_i=a_iw_1+b_iw_2$. Substituting this for $z_1{}^2+\cdots+z_N^{2}=0$, we get
$$
\Big(\sum_{i=1}^Na_i^2\Big)w_1^2+2\Big(\sum_{i=1}^Na_ib_i\Big)w_1w_2+\Big(\sum_{i=1}^Nb_i^2\Big)w_2^2=0
$$
so that
$$
\frac{w_2}{w_1}=\frac{-\sum_{i=1}^Na_ib_i\pm \sqrt{\Big(\sum_{i=1}^Na_ib_i\Big)^2-\Big(\sum_{i=1}^Na_i^2\Big)\Big(\sum_{i=1}^Nb_i^2\Big)}}{\sum_{i=1}^Nb_i^2}.
$$
If we denote by $D$ the square free part of the positive integer  
$$\Big(\sum_{i=1}^Na_i^2\Big)\Big(\sum_{i=1}^Nb_i^2\Big)-
\Big(\sum_{i=1}^Na_ib_i\Big)^2,
$$ 
then $x:=w_2/w_1\in \mathbb{Q}(\sqrt{-D})$. Thus (if $z_1\neq 0$),
$$
\frac{z_i}{z_1}=\frac{a_i+b_ix}{a_1+b_1x}\in \mathbb{Q}(\sqrt{-D}).
$$
Therefore 
$$
[z_1,\ldots,z_N]=[1,z_2/z_1,\ldots, z_N/z_1]\in P^{N-1}\big(\mathbb{Q}(\sqrt{-D})\big).
$$

Conversely suppose that $[z_1,\ldots,z_N]\in P^{N-1}\big(\mathbb{Q}(\sqrt{-D})\big)$. One may assume $z_i\in \mathbb{Q}(\sqrt{-D})$ without loss of generality. One may also assume $z_1,z_2$ are linearly independent over $\mathbb{R}$; otherwise every $z_i$ is a real scalar multiple of some $z$; say, $z_i=c_iz,~ c_i\in \mathbb{R}$, and hence $0=(c_1^2+\cdots+c_N^2)|z|^2$; thereby leading to a contradiction. Obviously there are rational numbers $\alpha_i,\beta_i$ such that $z_i=\alpha_iz_1+\beta_iz_2 $ ~$(i=3,\ldots, N)$. This implies that $\{z_1,\ldots,z_N\}$ generates a lattice in $\mathbb{C}$ (which is commensurable to the lattice generated by $z_1,z_2$. \hfill $\Box$

\medskip
 

\begin{thm} {\rm (1)} 
Let ${\bf z}\in C^1(X_0,\mathbb{C})$ be the building cochain of the standard realization of a 2D topological crystal over $X_0$. 
Then $[z_1, \ldots,z_N]\in {\bf Q}(X_0)\cap P^{N-1}\big(\mathbb{Q}(\sqrt{-D})\big)$. Conversely, for $[z_1, \ldots,z_N]\in {\bf Q}(X_0)\cap P^{N-1}\big(\mathbb{Q}(\sqrt{-D})\big)$, put 
$$
{\bf z}(e)=\begin{cases}
z_i & (e=e_i)\\
-z_i & (e=\overline{e_i}).
\end{cases}
$$
Then ${\bf z}$ is the building cochain of the standard realization of a 2D topological crystal. 

{\rm (2)} 
The set 
$$
{\bf Q}(X_0)\cap \bigcup_{D}P^{N-1}\big(\mathbb{Q}(\sqrt{-D})\big)
$$ 
is identified with the family of all similarity classes of standard realizations of 2-dimensional topological crystals over $X_0$. 

\end{thm}

It suffices to prove that for 
$$
[z_1,\ldots,z_N]\in {\bf Q}(X_0)\cap P^{N-1}\big(\mathbb{Q}(\sqrt{-D})\big),
$$
the rank of the image of $[{\bf z}]:H_1(X_0,\mathbb{Z})\longrightarrow \mathbb{C}$ is equal to two (note that ${\rm Image}~\![{\bf z}]$ is a subgroup of the lattice generated by $\{z_1,\ldots,z_N\}$). If the rank is one, then $\{z_1,\ldots,z_N\}$ gives rise to 1-dimensional (harmonic) realization, and hence every $z_i$ must be a real scalar multiple of some $z$; thereby a contradiction.
\hfill $\Box$

\medskip


We have interest in the description of $D$ in terms of $X_0$ and $H$.
 Consider the standard realization associated with 
$$
[z_1,\ldots,z_N]\in {\bf Q}(X_0)\cap P^{N-1}\big(\mathbb{Q}(\sqrt{-D})\big).
$$
One may assume $z_i\in \mathbb{Q}(\sqrt{-D})$. 
Let $T=\mathbb{Z}w_1+\mathbb{Z}w_2$ be the period lattice. Writing $w_i=a_i+b_i\sqrt{-D}\in \mathbb{Q}(\sqrt{-D})$, we observe that the energy of the standard realization associated with $[z_1,\ldots,z_N]$ is computed as
\begin{eqnarray*}
&&2{\rm vol}(\mathbb{C}/T)^{-1}(|z_1|^2+\cdots +|z_N|^2)=\frac{4}{|w_1\overline{w_2}-w_2\overline{w_1}|}(|z_1|^2+\cdots +|z_N|^2)\\
&&= \frac{2}{|a_1b_2-a_2b_1|\sqrt{D}}(|z_1|^2+\cdots +|z_N|^2).
\end{eqnarray*}
On the other hand, in view of (\ref{eq:energymini}), this is equal to
$$
4\kappa(X_0)^{-1/2}I(H)^{1/2}.
$$
Since $|z_1|^2+\cdots +|z_N|^2\in \mathbb{Q}$, we have 
\begin{prop}
$D$ is equal to the square free part of $\kappa(X_0)I(H)$.
\end{prop}

The following proposition is shown in a routine way.

\begin{prop}
If ${\bf Q}(X_0)$ has a $\mathbb{Q}(\sqrt{-D})$-rational point, then ${\bf Q}(X_0)\cap P^{N-1}\big(\mathbb{Q}(\sqrt{-D})\big)$ is dense in ${\bf Q}(X_0)$.
\end{prop}


\section{Examples}\label{sec:example}
We shall illustrate several examples.

\medskip

(1) The square lattice and the honeycomb have a special feature in the sense that they correspond to rational points of ``0-dimensional" quadrics; say:
 
In the case of the square lattice 
$$
\{[z_1,z_2]\in P^1(\mathbb{C});~z_1{}^2+z_2{}^2=0\} =[1,\pm\sqrt{-1}].
$$ 

In the case of the honeycomb 
\begin{eqnarray*}
&&\{[z_1,z_2,z_3]\in P^2(\mathbb{C});~z_1{}^2+z_2{}^2+z_3{}^2=0,~z_1+z_2+z_3=0\}\\
&&=\Big[1, \frac{-1\pm\sqrt{-3}}{2}, \frac{-1\mp\sqrt{-3}}{2}\Big].
\end{eqnarray*} 

\medskip

\begin{figure}[htbp]
\begin{center}
\includegraphics[width=.3\linewidth]
{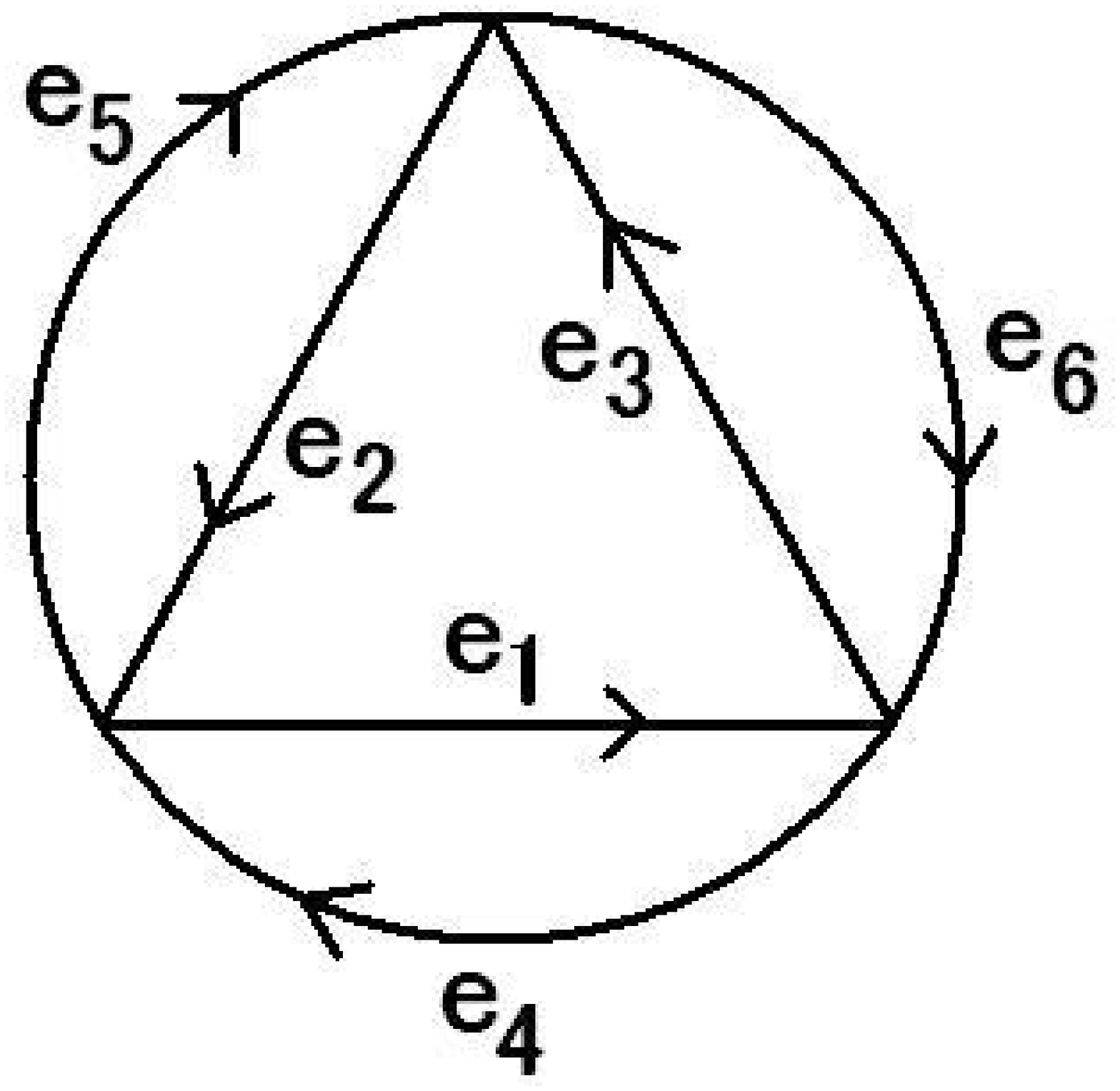}
\end{center}
\caption{Example}\label{fig:kagome}
\end{figure}

(2) The regular kagome lattice, which is, as an abstract graph,  a topological crystal over the graph depicted in Fig.~\!\ref{fig:kagome}, corresponds to the $\mathbb{Q}(\sqrt{-3})$-rational points 
$$
\Big[\frac{1\pm\sqrt{-3})}{2}, \frac{1\mp\sqrt{-3}}{2}, -1, \frac{1\pm\sqrt{-3}}{2}, \frac{1\mp\sqrt{-3}}{2}, -1\Big]
$$ 
of the 2-dimensional projective quadric defined over $\mathbb{Q}$
\begin{eqnarray*}
&&\{[z_1,z_2,z_3, z_4,z_5,z_6]\in P^5(\mathbb{C});~z_1{}^2+\cdots+z_6{}^2=0,\\
&& \qquad z_1+z_6=z_3+z_4,~ z_2+z_4=z_1+z_5,~z_3+z_5=z_2+z_6\}
\end{eqnarray*}
(the vanishing subgroup is $H=\mathbb{Z}(e_1+e_2+e_3)+\mathbb{Z}(e_4+e_5+e_6)$).

\medskip

\begin{figure}[htbp]
\begin{center}
\includegraphics[width=.6\linewidth]
{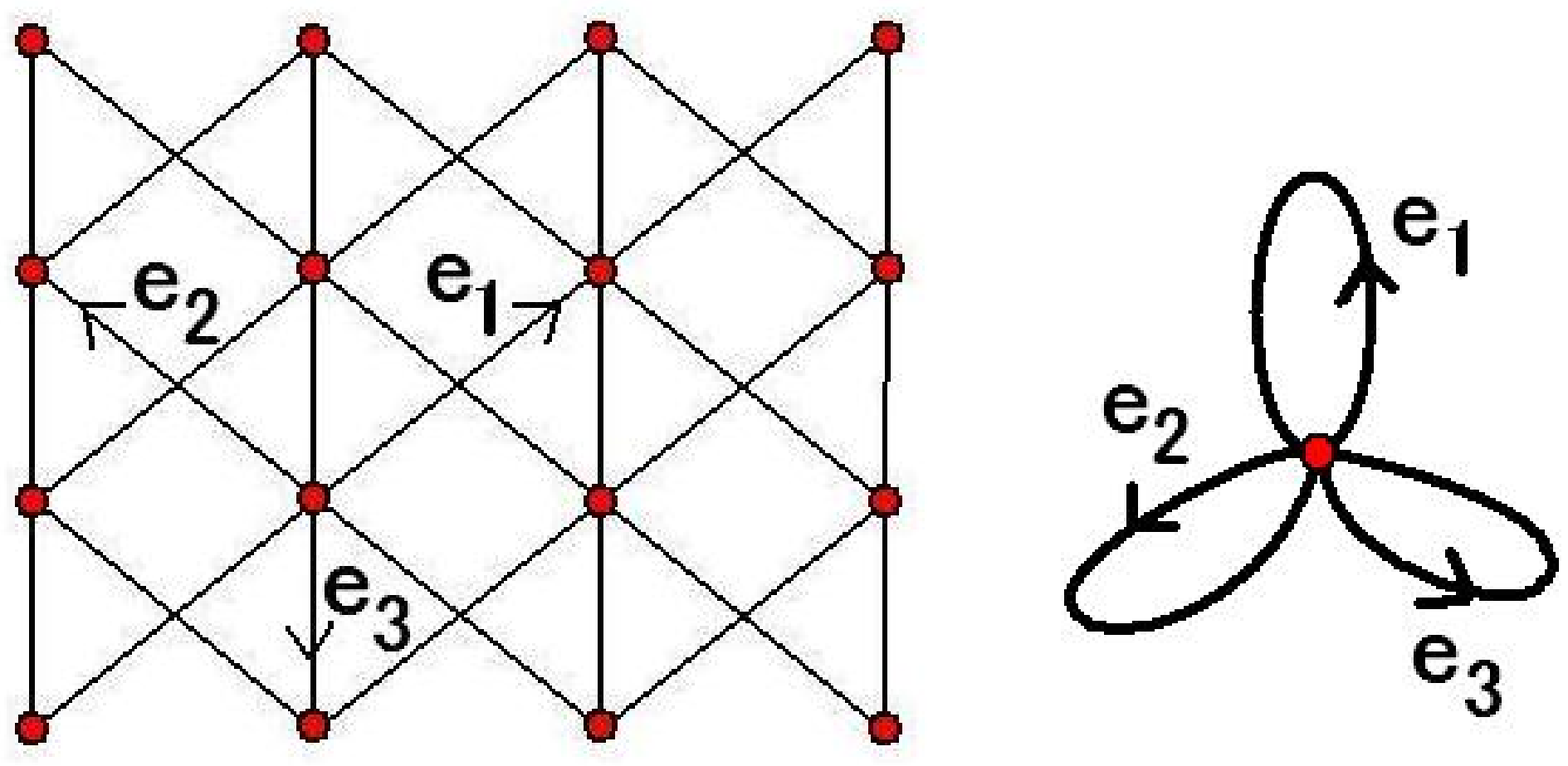}
\end{center}
\caption{Example}\label{fig:example1}
\end{figure}

(3) The equilateral triangular lattice, which is, as an abstract graph, a topological crystal over the 3-bouquet graph $B_3$,  corresponds to the $\mathbb{Q}(\sqrt{-3})$-rational points 
$$\Big[1, \frac{-1\pm\sqrt{-3}}{2}, \frac{-1\mp\sqrt{-3}}{2}\Big]$$ 
of the 
quadric 
$$\{[z_1,z_2,z_3]\in P^2(\mathbb{C});~z_1{}^2+z_2{}^2+z_3{}^2=0\}$$
(the vanishing group is $\mathbb{Z}(e_1+e_2+e_3)$; the orientation being illustrated in Fig.~\!\ref{fig:example1}).

\medskip

(4) The 2D crystal lattice in Fig.~\!\ref{fig:example1}, which is also a topological crystal over $B_3$ as an abstract graph, corresponds to $\mathbb{Q}(\sqrt{-6})$-rational points
$$
[3\pm\sqrt{-6}, -3\pm\sqrt{-6}, \mp\sqrt{-6}]
$$
of the quadric 
$$\{[z_1,z_2,z_3]\in P^2(\mathbb{C})|~z_1{}^2+z_2{}^2+z_3{}^2=0\}
$$
(the vanishing groups is $H=\mathbb{Z}(e_1+e_2+2e_3)$).

\medskip

\begin{figure}[htbp]
\begin{center}
\includegraphics[width=.7\linewidth]
{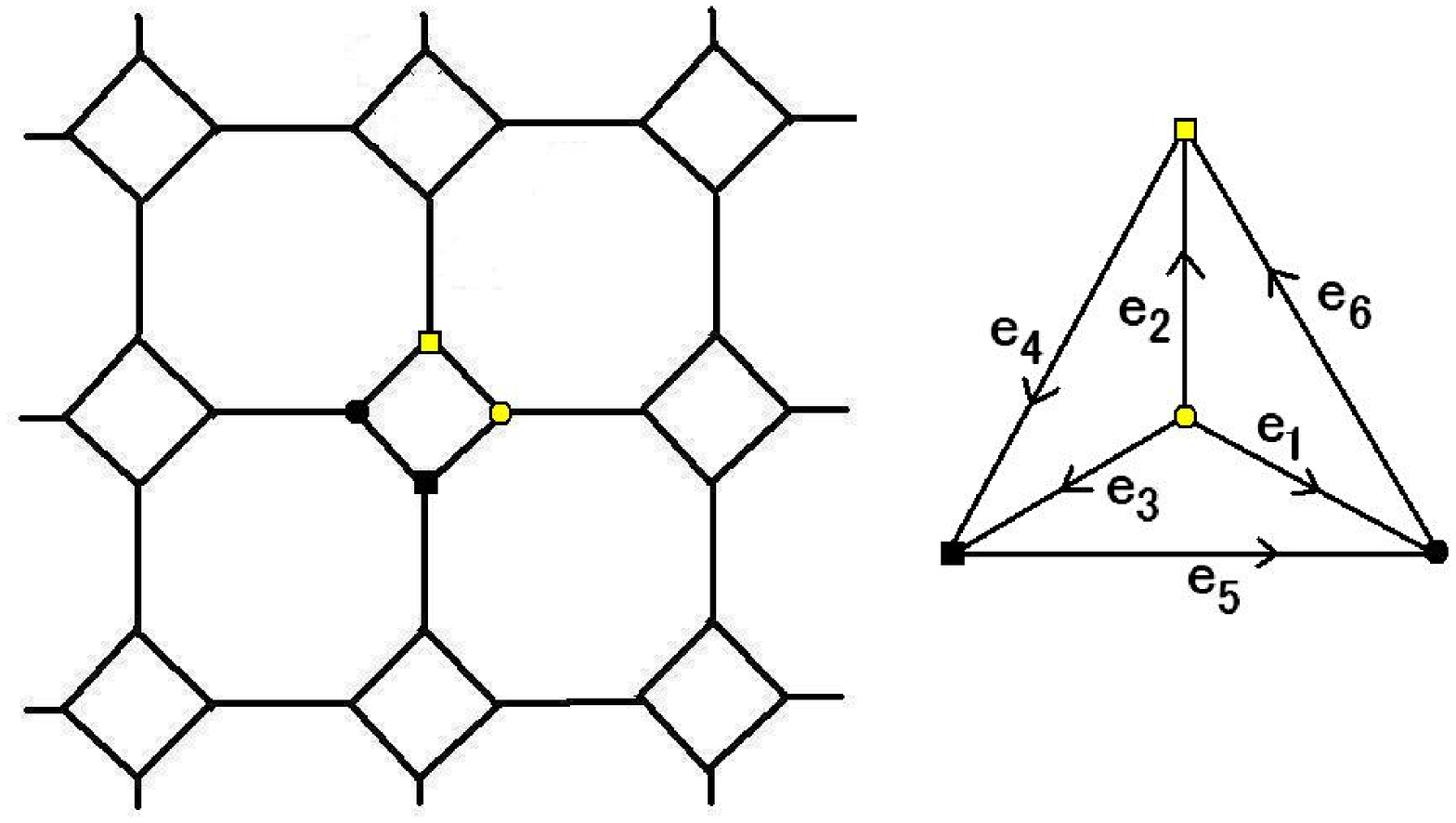}
\end{center}
\caption{Example}\label{fig:8-4}
\end{figure}

(5) The 2D crystal lattice in Fig.~\!\ref{fig:8-4} corresponds to $\mathbb{Q}(\sqrt{-1})$-rational points
$$
\Big[
1, \frac{-1\pm\sqrt{-1}}{2}, -\frac{1\pm\sqrt{-1}}{2}, \pm\sqrt{-1}, \frac{-1\pm\sqrt{-1}}{2}, \frac{1\pm\sqrt{-1}}{2}
\Big]
$$
of
the quadric
\begin{eqnarray*}
&&\{[z_1,\ldots,z_6]\in P^5(\mathbb{C});~z_1{}^2+\cdots+z_6{}^2=0,~z_1+z_2+z_3=0,\\
 &&\qquad\qquad z_3+z_4=z_5, z_1+z_5=z_6, z_2+z_6=z_4\}
\end{eqnarray*} 
(the vanishing group is $H=\mathbb{Z}(-e_2+e_3+e_5+e_6)$).

This crystalline pattern is observed when we look at the $K_4$ crystal (diamond twin) toward an suitable direction (see \cite{su2}, \cite{su5}).

\medskip

\begin{figure}[htbp]
\begin{center}
\includegraphics[width=.8\linewidth]
{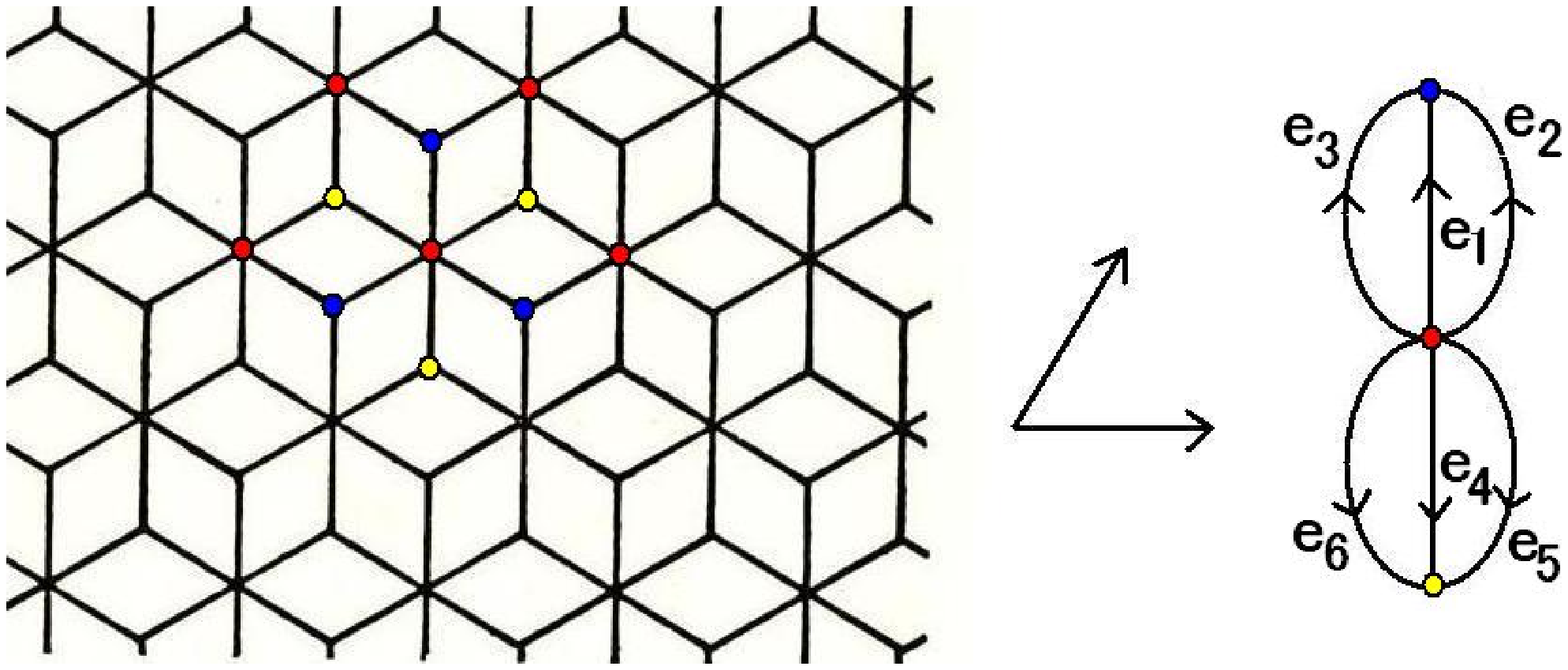}
\end{center}
\caption{Example}\label{fig:dice}
\end{figure}

(6) Figure \ref{fig:dice} is the so-called {\it dice lattice}. This corresponds to $\mathbb{Q}(\sqrt{-3})$-rational points
$$
\Big[1, \frac{-1\pm\sqrt{-3}}{2}, \frac{-1\mp\sqrt{-3}}{2}, -1, \frac{1\pm\sqrt{-3}}{2}, \frac{1\mp\sqrt{-3}}{2}\Big]
$$
of the quadric 
$$\{[z_1,\ldots,z_6]\in P^5(\mathbb{C})|~z_1{}^2+\cdots+z_6{}^2=0,~
z_1+z_2+z_3=0,~z_4+z_5+z_6=0\}
$$
(the vanishing groups is $H=\mathbb{Z}(e_1+\overline{e_3}+e_4+\overline{e_5})
+\mathbb{Z}(e_5+\overline{e_6}+e_3+\overline{e_2})$).

\begin{figure}[htbp]
\begin{center}
\includegraphics[width=.65\linewidth]
{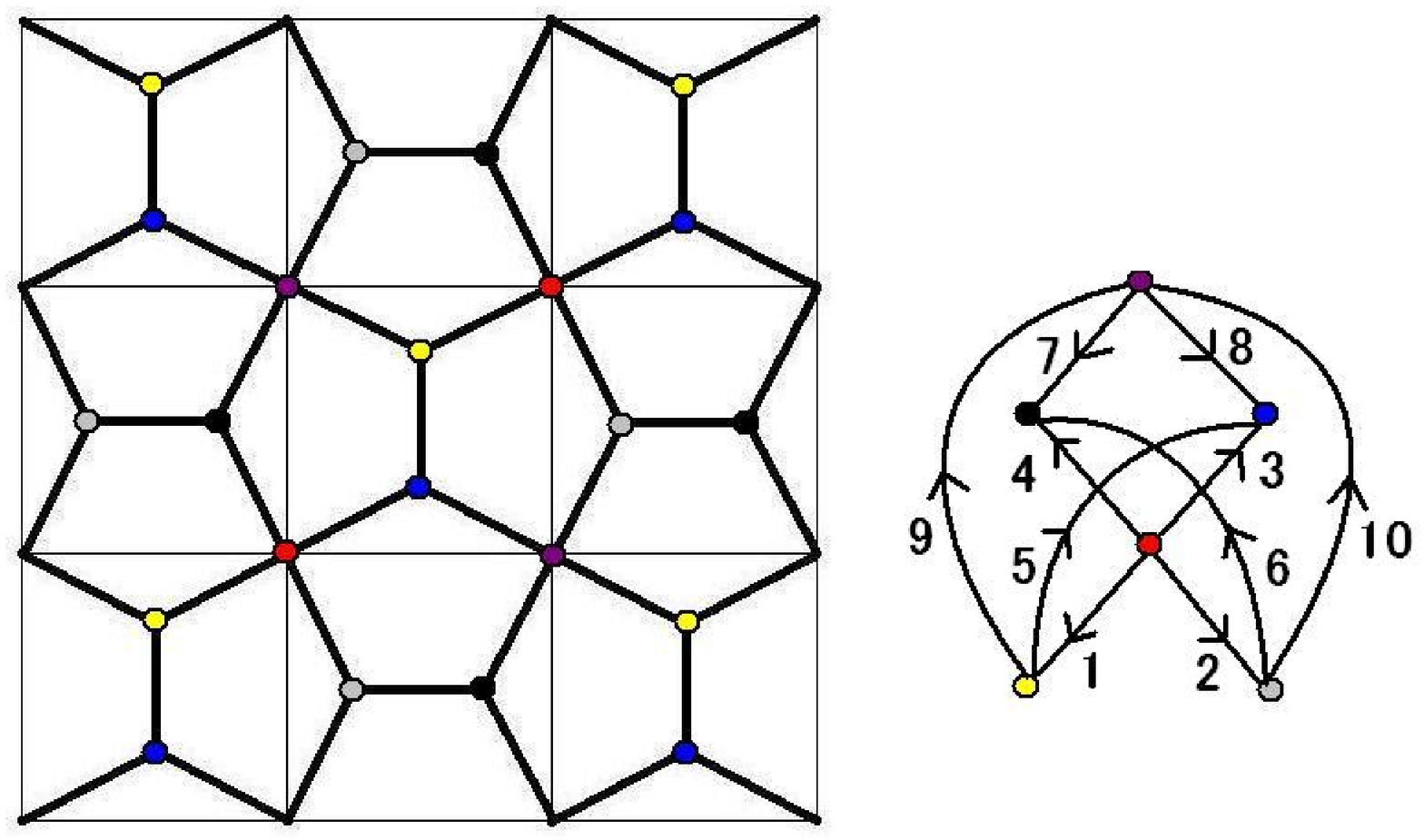}
\end{center}
\caption{Example}\label{fig:cairo1}
\end{figure}

(7) Figure \ref{fig:cairo1} is a tiling of pentagons with picturesque properties that has become known as the {\it Cairo pentagon}\footnote{It is given its name because several streets in Cairo are paved in this design (strictly speaking, it is a bit distorted). 
This is also called Macmahon's net, {\bf mcm}, and a 4-fold pentille (J. H. Conway).}. Its 1-skeleton is the standard realization of a topological crystal over the finite graph drawn on the right, which corresponds to $\mathbb{Q}(\sqrt{-1}))$-rational points
\begin{eqnarray*}
&&\Big[
-1\mp\frac{\sqrt{-1}}{2},
\frac{1}{2}\mp\sqrt{-1},
1\pm\frac{\sqrt{-1}}{2},
-\frac{1}{2}\pm\sqrt{-1},
\mp\sqrt{-1},\\
&&\qquad 1,
-\frac{1}{2}\mp\sqrt{-1},
-1\pm\frac{\sqrt{-1}}{2},
-1\pm\frac{\sqrt{-1}}{2},
-\frac{1}{2}\mp\frac{\sqrt{-1}}{2}
\Big]
\end{eqnarray*} 
of the quadric
\begin{eqnarray*}
&&\{[z_1,\ldots,z_{10}]\in P^{9}(\mathbb{C})|~ z_1{}^2+\cdots+z_{10}{}^2=0,~
 z_1=z_5+z_9,~z_2=z_6+z_{10},\\
&&\qquad\qquad \qquad z_1+z_2+z_3+z_4=0,~z_4+z_6+z_7=0,~ z_3+z_5+z_8=0,\\
&&\qquad\qquad \qquad z_9+z_{10}=z_7+z_8
\}.
\end{eqnarray*}

\section{Tilings}
In general, a {\it periodic tiling}, symbolically written as $(T,L)$, is a {\it tessellation} of figures (tiles) in the plane $\mathbb{R}^2$ which is periodic with respect to the translational action by a lattice group $L$. Two tiles $D$ and $D'$ in $(T,L)$ are said to be {\it equivalent} if $D'=D+\sigma$ for some $\sigma\in L$. We denote by $f_{T,L}$ the number of equivalence classes of tiles, and let $D_1,\ldots,D_{f_{T,L}}$ be representatives of equivalence classes, which we call {\it fundamental tiles}\footnote{Some of fundamental tiles can be congruent. For instance, fundamental tiles of the equilateral lattice with respect to the maximal periodic lattice are two congruent equilateral triangles.} of $(T,L)$. 

\begin{figure}[htbp]
\begin{center}
\includegraphics[width=.5\linewidth]
{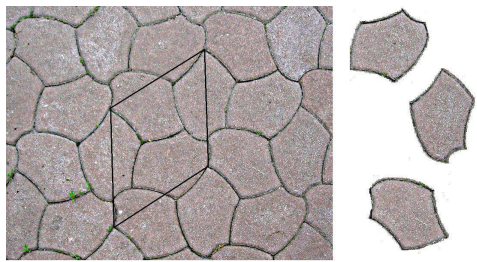}
\end{center}
\caption{A tiling and its fundamental tiles}\label{fig:fund}
\end{figure}

Figure \ref{fig:fund} is the picture of the tiling we find on a street pavement in Zakopane, Poland\footnote{From Wikipedia.}. 
This tiling has three fundamental tiles depicted on the right.

As for topological shapes of tilings, we have 

\begin{thm}\label{thm1} For any natural number $f$, there are only finitely many homeomorphic classes of tilings $(T,L)$ such that $f_{T,L}=f$. Here we say that two tilings $(T_1,L_1)$ and $(T_2,L_2)$ are {\it homemorphic} if there is a homemorphism $\varphi:\mathbb{R}^2\longrightarrow \mathbb{R}^2$ and an isomorphism $\psi:L_1\longrightarrow L_2$ satisfying 

(1) $\varphi (x+\sigma)=\varphi(x)+\psi (\sigma), \quad (x\in \mathbb{R}^2,~\sigma \in L_1)$,

(2) $\varphi(T_1)=T_2$.

\end{thm}

The following theorem is deduced from the proof of this theorem.

\begin{thm}\label{thm2} There are only finitely many 2D topological crystals over $X_0$ whose standard realizations yield tilings. In other words, there are only finitely many rational points in ${\bf Q}(X_0)$ which correspond to tilings.
\end{thm}

Let $X$ be the 1-skeleton of a given periodic tiling $(T,L)$, and consider the quotient graph $X_0=X/L$. Then $X$ as an abstract graph is a 2D topological crystal over $X_0$ (the net $X$ is not necessarily a crystal net since the edges are allowed to be curved). We shall say that $(T,L)$ is a {\it tiling with the base graph} $X_0$. Denote by $\omega:X\longrightarrow X_0$ the covering map. Note that the degree of any vertex (node) $x$ in $X_0$ is greater than or equal to 3. 

The covering map $\pi$ of $\mathbb{R}^2$ onto the the 2-dimensional torus $\mathbb{R}^2/L$ induces a cellular decomposition of $\mathbb{R}^2/L$. The 1-skeleton of this cellular decomposition is just the finite graph $X_0$ realized in $\mathbb{R}^2/L$ (conversely, a cellular decomposition of the torus yields a tiling).

The proof of Theorem \ref{thm1} is divided into 4 steps each of which is fairly elementary. 

\medskip

(1) The number of 2-cells in $\mathbb{R}^2/L$ is equal to $f=f_{T,L}$, which coincides with $v-e~(=b_1(X_0)-1)$, where $v$ (resp. $e$) is the number of vertices (resp. edges) in $X_0$. This is because $v-e+f$ is the Euler number of the torus so that $v-e+f=0$.

\smallskip

(2) Let $D_1,\ldots,D_{b-1}$ $(b=b_1(X_0))$ be fundamental tiles. Then $\pi(D_1),\ldots,$ $\pi(D_{b-1})$ are all 2-cells in the torus. Topologically $D_i$ is identified with a polygon with $k_i$ sides where $k_i$ is the number of vertices on $D_i$. If we put $N_{T,L}=\max\{k_1,\ldots,k_{b-1}\}$, then $N_{T,L}\leq 2e$.

To show this, let $D_i$ is a tile having $N_{T,L}$ edges on $\partial D_i$. 
If $N_{T,L}>2e$, then there exist at least three edges in $\partial D_i$ which are mapped to an edge in $X_0$ by $\omega$. This should not happen since the torus is non-singular.  Therefore we get the claim.

\smallskip

(3) For a fixed integer $b\geq 0$, there are only finitely many finite graphs $X_0$ $($up to isomorphisms$)$ such that ${\rm deg}~\hspace{-0.04cm}x\geq 3$ for all vertex $x$ and $b_1(X_0)=b$.

Indeed, 
$$
3v\leq \sum_{x}{\rm deg}~\!x=2e,
$$
so $3v\leq 2e$ and $v\leq 2(b-1)$ (use $v-e=1-b$), and also $e\leq 3b-3$. Therefore the number of vertices and edges is bounded.

\smallskip

(4) Given a finite graph $X_0$ and polygonal 2-cells $D_i$ with $k_i$ sides ($i=1,\ldots,f$, there are only finitely many ways to attach $D_i$'s to $X_0$ to make a cell complex whose underlying space is the torus (more generally a closed surface). 

\medskip
Putting altogether, we complete the proof of 
Theorem \ref{thm1}. Theorem \ref{thm2} is also a consequence of (1), (2), (4) since if a 2D topological crystal $X$ over $X_0$ yields a tiling $(T,L)$, then $f_{T,L}=b_1(X_0)-1$ and $N_{T,L}\leq 6\big(b_1(X_0)-1\big)$.

\medskip

We shall go a bit further. Let $H$ be the vanishing group for the 2D topological crystal $X$ over $X_0$ which corresponds to a tiling $(T,L)$.

\begin{lemma}
If we put $c_i=\partial D_i$ and give the counter-clockwise orientation on $c_i$, then $\omega(c_{b-1})=-(\omega(c_1)+\cdots+\omega(c_{b-2}))$, and 
$\omega(c_1),\ldots,\omega(c_{b-2})$ form a $\mathbb{Z}$-basis of $H$. 

\end{lemma}

To prove the last assertion, let $c$ be a closed path in $X_0$ whose homology class is in $H$. The (any) lifting $\widehat{c}$ of $c$ in $X$ is also closed. As a 1-chain, we may write 
$$
\widehat{c}=\sum_{i=1}^{b-1}n_i\partial D_i,
$$ 
where $n_i$ denotes the winding number of $\widehat{c}$ around a interior point of $D_i$. Thus 
$$
c=\omega(\widehat{c})=\sum_{i=1}^{b-1}n_i\omega(c_i),
$$ 
which implies that $H$ is generated by $\omega(c_1),\ldots,\omega(c_{b-2})$.

Next suppose that 
$$
\sum_{i=1}^{b-2}m_i\omega(c_i)=0\qquad (m_i\in \mathbb{Z}).
$$
Expressing the left hand side as a 1-chain, and taking a look at the coefficients of directed edges in $X_0$, we find that if  
$\omega(c_i)$ and $\omega(c_j)$ share an edge, then $m_i=m_j$, and that if $\omega(c_i)$ and $\omega(c_{b-1})$ share an edge, then $m_i=0$, from which it follows that $m_1=\cdots=m_{b-2}=0$. Thus we conclude that $\omega(c_1),\ldots,\omega(c_{b-2})$ comprise a $\mathbb{Z}$-basis of $H$.
\smallskip

It is an interesting problem to list all homeomorphic classes of tilings $(T,L)$ with $f_{T,L}=f$. To handle this problem (not yet having been solved completely), it may be useful to introduce the notion of {\it height} of vanishing groups as follows\footnote{This notion is introduced for topological crystal of arbitrary dimension (\cite{su5}).}. 

Take an orientation of $X_0$ and define the norm $\|\alpha\|_1$ of a $1$-chain\index{$1$-chain} $\alpha=\displaystyle\sum_{e\in E_0^o}a_ee$ by setting
$$
\|\alpha\|_1=\sum_{e\in E_0^o}|a_e|,
$$
where $E_0^o$ is the set of directed edges for the orientation (it should be noted that $\|\alpha\|_1$ does not depend on the choice of an orientation). For a subset $S=\{\alpha_{1},\ldots, \alpha_{b-2}\}$ of $H_1(X_0,\mathbb{Z})$ which forms a $\mathbb{Z}$-basis of a vanishing subgroup, we put
$$
h(S)=\max (\|\alpha_{1}\|_1,\ldots, \|\alpha_{b-2}\|_1),
$$ 
and define the {\it height}\index{height} of $H$ by 
$
h(H)=\min_{S}h(S),
$
where $S$ runs over all $\mathbb{Z}$-bases of $H$. 
Consider two sets
$$ 
R_1=\{S|~h(S)\leq h\},\quad R_2=\{H|~h(H)\leq h\}.
$$
Certainly $R_1$ is a finite set. The correspondence 
$
S\mapsto H~(\text{generated by S})
$
yields a surjective map of $R_1$ onto $R_2$. Thus we get: 

\begin{thm}\label{thm:finiteness} There are only a finite number of vanishing subgroups $H$ of the homology group\index{homology group} $H_1(X_0,\mathbb{Z})$ such that 

\smallskip

$(1)$ ${\rm rank}~\hspace{-0.05cm}H_1(X_0,\mathbb{Z})/H=2$,

\smallskip

$(2)$ $h(H)\leq h$.

\end{thm} 

Suppose now that $X_0$ is the base graph for a tiling $(T,L)$, and let $H$ be the vanishing group corresponding to the topological crystal $X$ associated with $(T,L)$. Then, using the notations above, we have
$$
\|\omega(c_i)\|_1\leq k_i\leq N_{T,L}
$$
so that $h(H)\leq N_{T,L}\leq 6\big(b_1(X_0)-1\big)$ (this gives another proof of Theorem \ref{thm2}). Thus in order to list homeomorphic classes of tiling $(T,L)$ with $f_{T,L}=f$, we first enumerate finite graphs $X_0$ with $b_1(X_0)=f+1$, and then for such $X_0$, we enumerate vanishing groups $H\subset H_1(X_0,\mathbb{Z})$ such that $h(H)\leq 6f$ and check whether $X=X_{0}^{\rm ab}/H$ gives a tiling or not.


\end{document}